\documentclass{amsart}
\usepackage{amsmath}
\usepackage{amssymb}
\usepackage{amsfonts}
\usepackage{cite}

\usepackage{geometry} 

\usepackage[table,usenames,dvipsnames]{xcolor}

\usepackage{gnuplot-lua-tikz}
\usepackage{tikz}
\usetikzlibrary{matrix,arrows,positioning,automata}
\pgfdeclarelayer{background layer}
\pgfsetlayers{background layer,main}

\definecolor{darkblue}{rgb}{0.0,0.0,0.3}
\definecolor{gray9}{gray}{0.9}
\definecolor{gray8}{gray}{0.8}
\definecolor{gray7}{gray}{0.7}
\definecolor{gray6}{gray}{0.6}
\definecolor{grey}{gray}{0.8}

\newcommand{\myrowcolour}{\rowcolor[gray]{0.925}}

\newcommand{\BinRel}{\text{B}}
\newcommand{\InvMon}{\mathcal I}
\newcommand{\DualInvMon}{{\mathcal I}^*}
\newcommand{\FullTrans}{\mathcal T}
\newcommand{\Symmetric}{\mathcal S}
\newcommand{\Brauer}{\mathfrak B}
\newcommand{\Partition}{\mathcal P}
\newcommand{\PartBinRel}{\Partition\BinRel}
\newcommand{\TemperleyLieb}{\text{TL}}
\newcommand{\PartialTrans}{\text{P}\FullTrans}
\newcommand{\Sub}{\mathbf{Sub}}

\newcommand{\cD}{\mathcal D}
\newcommand{\cR}{\mathcal R}
\newcommand{\cL}{\mathcal L}

\newcommand{\cJ}{\mathcal J}

\newcommand{\id}{1}

\newcommand{\GAP}{\textsc{Gap}}

\newcommand{\semigroups}{\textsc{Semigroups}}
\newcommand{\subsemi}{\textsc{SubSemi}}
\newcommand{\smallsemi}{\textsc{SmallSemi}}

\newcommand{\todo}[1]{}
\newcommand{\gap}{\vskip6pt}

\DeclareMathOperator{\Dom}{Dom}
\DeclareMathOperator{\Codom}{Codom}

\begin{document}

\title{Finite Diagram Semigroups: Extending the Computational Horizon}
\author[J. East, A. Egri-Nagy, A.R. Francis, J. D. Mitchell]{James East$^1$, Attila Egri-Nagy$^{1}$, Andrew R. Francis$^1$, James D. Mitchell$^2$}
\address{$^1$Centre for Research in Mathematics, School of Computing, Engineering and Mathematics, University of Western Sydney (Parramatta Campus), Locked Bag 1797, Penrith, NSW 2751, Australia}
\address{$^2$ Mathematical Institute, University of St Andrews, North Haugh, St Andrews, Fife, KY16 9SS, Scotland}
\email{J.East@uws.edu.au,\ A.Egri-Nagy@uws.edu.au,\ A.Francis@uws.edu.au,\ jdm3@st-and.ac.uk}

\maketitle

\begin{abstract}

Diagram semigroups are interesting algebraic and combinatorial objects, several types of them originating from questions in computer science and in physics.
Here we describe diagram semigroups in a general framework and extend our computational knowledge of them.
The generated data set is replete with surprising observations raising many open questions for further theoretical research.  
\end{abstract}

\section{Introduction}

For studying finite structures  it is helpful to generate small examples by computer programs.
By investigating these sample objects we can formulate new hypotheses and falsify conjectures by counterexamples.
At a given time, the available computing power and the state of the art algorithms define a limit on the size of the examples we can investigate.
This limit we call the \emph{computational horizon}, similar to the cosmological horizon determined by the size of our observable physical universe.
 The underlying assumption in both fields is that we can see enough within our limits to enable us to go beyond by theoretical reasoning, i.e.~to have enough observational data to construct valid theories. Here we aim to extend the database of small degree diagram semigroups.

\section{Diagram Semigroups}

Diagram representations of finite semigroups are described by fundamental mathematical objects such as relations and functions, and so they often arise naturally in mathematical theories.
The original interest came from algebras with a basis whose elements can be multiplied diagrammatically (e.g.~\cite{Brauer1937,Jones1994,Martin1994}).

In what follows, we define diagram semigroups in a logical order (as opposed to a historical order) starting with the most general diagram type and define each type by a set of constraints.
The conceptual origin of diagram semigroups is the notion of binary relation on a set $A$.
Such a relation, a subset of $A\times A$, can be represented as a graph by a set of directed edges between the elements of $A$ as vertices.
In order to make the graph of a binary relation into a diagram that can be combined with other diagrams we partition $A$ into two parts.
The \emph{domain} $\Dom(A)$ and \emph{codomain} $\Codom(A)$ are the ``interfaces'' for combining diagrams.
In general we can talk about $(n,m)$ diagrams, where $n$ is the size of the domain and $m$ is the size of the codomain and $|A|=n+m$. 
Here, we restrict our attention to $(n,n)$ diagrams, so we partition $A$ into to equal sized parts.

\subsection{Partitioned Binary Relations}

Partitioned binary relations are the most general type of diagrams we consider, although historically it was the last to be defined. For the formal definitions and its categorical context see \cite{PartBinRel2013}.

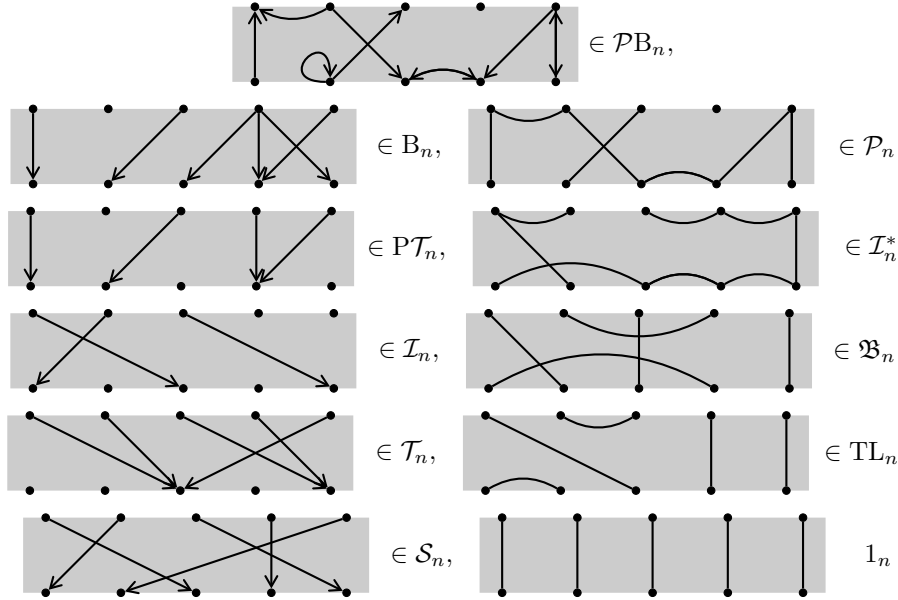
\begin{figure}[h]
\begin{center}
\begin{tikzpicture}
\tikzstyle{blackdot}=[draw=black,circle,fill=black,inner sep=1pt]
\tikzstyle{arrow}=[thick,->,>=angle 60]
\node [blackdot] at (0,0) (u1) {};
\node [blackdot,right of=u1] (u2) {};
\node [blackdot,right of=u2] (u3) {};
\node [blackdot,right of=u3] (u4) {};
\node [blackdot,right of=u4] (u5) {};
\node [blackdot,below of=u1] (d1) {};
\node [blackdot,below of=u2] (d2) {};
\node [blackdot,below of=u3] (d3) {};
\node [blackdot,below of=u4] (d4) {};
\node [blackdot,below of=u5] (d5) {};
\node at (5,-.5) {$\in\Partition\BinRel_n$,};
\draw [arrow] (d1) edge (u1);
\draw [arrow,bend left] (u2) edge (u1);
\draw [arrow] (u2) edge (d3);
\draw [arrow] (d2) edge (u3);
\draw [arrow] (d2) to [out=180,in=90,looseness=24] (d2);
\draw [arrow,bend left] (d3) edge (d4);
\draw [arrow,bend right] (d4) edge (d3);
\draw [arrow] (u5) edge (d5);
\draw [arrow] (d5) edge (u5);
\draw [arrow] (u5) edge (d4);
\begin{pgfonlayer}{background layer}
\fill  [grey] plot (-.3,0) rectangle (4.3,-1);
\end{pgfonlayer}
\end{tikzpicture}
\gap
\begin{tikzpicture}
\tikzstyle{blackdot}=[draw=black,circle,fill=black,inner sep=1pt]
\tikzstyle{arrow}=[thick,->,>=angle 60]
\node [blackdot] at (0,0) (u1) {};
\node [blackdot,right of=u1] (u2) {};
\node [blackdot,right of=u2] (u3) {};
\node [blackdot,right of=u3] (u4) {};
\node [blackdot,right of=u4] (u5) {};
\node [blackdot,below of=u1] (d1) {};
\node [blackdot,below of=u2] (d2) {};
\node [blackdot,below of=u3] (d3) {};
\node [blackdot,below of=u4] (d4) {};
\node [blackdot,below of=u5] (d5) {};
\node at (5,-.5) {$\in\BinRel_n$,};
\draw [arrow] (u1) edge (d1);
\draw [arrow] (u3) edge (d2);
\draw [arrow] (u4) edge (d4);
\draw [arrow] (u4) edge (d5);
\draw [arrow] (u4) edge (d3);
\draw [arrow] (u5) edge (d4);
\begin{pgfonlayer}{background layer}
\fill  [grey] plot (-.3,0) rectangle (4.3,-1);
\end{pgfonlayer}
\end{tikzpicture}
\begin{tikzpicture}
\tikzstyle{blackdot}=[draw=black,circle,fill=black,inner sep=1pt]
\tikzstyle{arrow}=[thick]
\node [blackdot] at (0,0) (u1) {};
\node [blackdot,right of=u1] (u2) {};
\node [blackdot,right of=u2] (u3) {};
\node [blackdot,right of=u3] (u4) {};
\node [blackdot,right of=u4] (u5) {};
\node [blackdot,below of=u1] (d1) {};
\node [blackdot,below of=u2] (d2) {};
\node [blackdot,below of=u3] (d3) {};
\node [blackdot,below of=u4] (d4) {};
\node [blackdot,below of=u5] (d5) {};
\node at (5,-.5) {$\in\Partition_n$};
\draw [arrow] (u1) edge (d1);
\draw [arrow,bend left] (u2) edge (u1);
\draw [arrow] (u2) edge (d3);
\draw [arrow] (d2) edge (u3);
\draw [arrow,bend left] (d3) edge (d4);
\draw [arrow,bend right] (d4) edge (d3);
\draw [arrow] (u5) edge (d5);
\draw [arrow] (d5) edge (u5);
\draw [arrow] (u5) edge (d4);
\begin{pgfonlayer}{background layer}
\fill  [grey] plot (-.3,0) rectangle (4.3,-1);
\end{pgfonlayer}
\end{tikzpicture}
\gap
\begin{tikzpicture}
\tikzstyle{blackdot}=[draw=black,circle,fill=black,inner sep=1pt]
\tikzstyle{arrow}=[thick,->,>=angle 60]
\node [blackdot] at (0,0) (u1) {};
\node [blackdot,right of=u1] (u2) {};
\node [blackdot,right of=u2] (u3) {};
\node [blackdot,right of=u3] (u4) {};
\node [blackdot,right of=u4] (u5) {};
\node [blackdot,below of=u1] (d1) {};
\node [blackdot,below of=u2] (d2) {};
\node [blackdot,below of=u3] (d3) {};
\node [blackdot,below of=u4] (d4) {};
\node [blackdot,below of=u5] (d5) {};
\node at (5,-.5) {$\in\PartialTrans_n$,};
\draw [arrow] (u1) edge (d1);
\draw [arrow] (u3) edge (d2);
\draw [arrow] (u4) edge (d4);
\draw [arrow] (u5) edge (d4);
\begin{pgfonlayer}{background layer}
\fill  [grey] plot (-.3,0) rectangle (4.3,-1);
\end{pgfonlayer}
\end{tikzpicture}
\begin{tikzpicture}
\tikzstyle{blackdot}=[draw=black,circle,fill=black,inner sep=1pt]
\tikzstyle{arrow}=[thick]
\node [blackdot] at (0,0) (u1) {};
\node [blackdot,right of=u1] (u2) {};
\node [blackdot,right of=u2] (u3) {};
\node [blackdot,right of=u3] (u4) {};
\node [blackdot,right of=u4] (u5) {};
\node [blackdot,below of=u1] (d1) {};
\node [blackdot,below of=u2] (d2) {};
\node [blackdot,below of=u3] (d3) {};
\node [blackdot,below of=u4] (d4) {};
\node [blackdot,below of=u5] (d5) {};
\node at (5,-.5) {$\in\DualInvMon_n$};
\draw [arrow] (u1) edge (d2);
\draw [arrow,bend left] (u2) edge (u1);
\draw [arrow,bend left] (u4) edge (u3);
\draw [arrow,bend left] (u5) edge (u4);
\draw [arrow,bend left] (d3) edge (d4);
\draw [arrow,bend right] (d3) edge (d1);
\draw [arrow,bend right] (d4) edge (d3);
\draw [arrow,bend right] (d5) edge (d4);
\draw [arrow] (u5) edge (d5);

\begin{pgfonlayer}{background layer}
\fill  [grey] plot (-.3,0) rectangle (4.3,-1);
\end{pgfonlayer}
\end{tikzpicture}
\gap
\begin{tikzpicture}
\tikzstyle{blackdot}=[draw=black,circle,fill=black,inner sep=1pt]
\tikzstyle{arrow}=[thick,->,>=angle 60]
\node [blackdot] at (0,0) (u1) {};
\node [blackdot,right of=u1] (u2) {};
\node [blackdot,right of=u2] (u3) {};
\node [blackdot,right of=u3] (u4) {};
\node [blackdot,right of=u4] (u5) {};
\node [blackdot,below of=u1] (d1) {};
\node [blackdot,below of=u2] (d2) {};
\node [blackdot,below of=u3] (d3) {};
\node [blackdot,below of=u4] (d4) {};
\node [blackdot,below of=u5] (d5) {};
\node at (5,-.5) {$\in\InvMon_n$,};
\draw [arrow] (u1) edge (d3);
\draw [arrow] (u2) edge (d1);
\draw [arrow] (u3) edge (d5);
\begin{pgfonlayer}{background layer}
\fill  [grey] plot (-.3,0) rectangle (4.3,-1);
\end{pgfonlayer}
\end{tikzpicture}
\begin{tikzpicture}
\tikzstyle{blackdot}=[draw=black,circle,fill=black,inner sep=1pt]
\tikzstyle{arrow}=[thick]
\node [blackdot] at (0,0) (u1) {};
\node [blackdot,right of=u1] (u2) {};
\node [blackdot,right of=u2] (u3) {};
\node [blackdot,right of=u3] (u4) {};
\node [blackdot,right of=u4] (u5) {};
\node [blackdot,below of=u1] (d1) {};
\node [blackdot,below of=u2] (d2) {};
\node [blackdot,below of=u3] (d3) {};
\node [blackdot,below of=u4] (d4) {};
\node [blackdot,below of=u5] (d5) {};
\node at (5,-.5) {$\in\Brauer_n$};
\draw [arrow] (u1) edge (d2);
\draw [arrow,bend right] (u2) edge (u4);
\draw [arrow,bend right] (d4) edge (d1);
\draw [arrow] (u3) edge (d3);
\draw [arrow] (u5) edge (d5);

\begin{pgfonlayer}{background layer}
\fill  [grey] plot (-.3,0) rectangle (4.3,-1);
\end{pgfonlayer}
\end{tikzpicture}
\gap 
\begin{tikzpicture}
\tikzstyle{blackdot}=[draw=black,circle,fill=black,inner sep=1pt]
\tikzstyle{arrow}=[thick,->,>=angle 60]
\node [blackdot] at (0,0) (u1) {};
\node [blackdot,right of=u1] (u2) {};
\node [blackdot,right of=u2] (u3) {};
\node [blackdot,right of=u3] (u4) {};
\node [blackdot,right of=u4] (u5) {};
\node [blackdot,below of=u1] (d1) {};
\node [blackdot,below of=u2] (d2) {};
\node [blackdot,below of=u3] (d3) {};
\node [blackdot,below of=u4] (d4) {};
\node [blackdot,below of=u5] (d5) {};
\node at (5,-.5) {$\in\FullTrans_n$,};
\draw [arrow] (u1) edge (d3);
\draw [arrow] (u2) edge (d3);
\draw [arrow] (u3) edge (d5);
\draw [arrow] (u4) edge (d5);
\draw [arrow] (u5) edge (d3);
\begin{pgfonlayer}{background layer}
\fill  [grey] plot (-.3,0) rectangle (4.3,-1);
\end{pgfonlayer}
\end{tikzpicture}
\begin{tikzpicture}
\tikzstyle{blackdot}=[draw=black,circle,fill=black,inner sep=1pt]
\tikzstyle{arrow}=[thick]
\node [blackdot] at (0,0) (u1) {};
\node [blackdot,right of=u1] (u2) {};
\node [blackdot,right of=u2] (u3) {};
\node [blackdot,right of=u3] (u4) {};
\node [blackdot,right of=u4] (u5) {};
\node [blackdot,below of=u1] (d1) {};
\node [blackdot,below of=u2] (d2) {};
\node [blackdot,below of=u3] (d3) {};
\node [blackdot,below of=u4] (d4) {};
\node [blackdot,below of=u5] (d5) {};
\node at (5,-.5) {$\in\TemperleyLieb_n$};
\draw [arrow] (u1) edge (d3);
\draw [arrow,bend right] (u2) edge (u3);
\draw [arrow,bend right] (d2) edge (d1);
\draw [arrow] (u4) edge (d4);
\draw [arrow] (u5) edge (d5);

\begin{pgfonlayer}{background layer}
\fill  [grey] plot (-.3,0) rectangle (4.3,-1);
\end{pgfonlayer}
\end{tikzpicture}
\gap
\begin{tikzpicture}
\tikzstyle{blackdot}=[draw=black,circle,fill=black,inner sep=1pt]
\tikzstyle{arrow}=[thick,->,>=angle 60]
\node [blackdot] at (0,0) (u1) {};
\node [blackdot,right of=u1] (u2) {};
\node [blackdot,right of=u2] (u3) {};
\node [blackdot,right of=u3] (u4) {};
\node [blackdot,right of=u4] (u5) {};
\node [blackdot,below of=u1] (d1) {};
\node [blackdot,below of=u2] (d2) {};
\node [blackdot,below of=u3] (d3) {};
\node [blackdot,below of=u4] (d4) {};
\node [blackdot,below of=u5] (d5) {};
\node at (5,-.5) {$\in\Symmetric_n$,};
\draw [arrow] (u1) edge (d3);
\draw [arrow] (u2) edge (d1);
\draw [arrow] (u3) edge (d5);
\draw [arrow] (u4) edge (d4);
\draw [arrow] (u5) edge (d2);
\begin{pgfonlayer}{background layer}
\fill  [grey] plot (-.3,0) rectangle (4.3,-1);
\end{pgfonlayer}
\end{tikzpicture}
\begin{tikzpicture}
\tikzstyle{blackdot}=[draw=black,circle,fill=black,inner sep=1pt]
\tikzstyle{arrow}=[thick]
\node [blackdot] at (0,0) (u1) {};
\node [blackdot,right of=u1] (u2) {};
\node [blackdot,right of=u2] (u3) {};
\node [blackdot,right of=u3] (u4) {};
\node [blackdot,right of=u4] (u5) {};
\node [blackdot,below of=u1] (d1) {};
\node [blackdot,below of=u2] (d2) {};
\node [blackdot,below of=u3] (d3) {};
\node [blackdot,below of=u4] (d4) {};
\node [blackdot,below of=u5] (d5) {};
\node at (5,-.5) {$1_n$};
\draw [arrow] (u1) edge (d1);
\draw [arrow] (u2) edge (d2);
\draw [arrow] (u3) edge (d3);
\draw [arrow] (u4) edge (d4);
\draw [arrow] (u5) edge (d5);

\begin{pgfonlayer}{background layer}
\fill  [grey] plot (-.3,0) rectangle (4.3,-1);
\end{pgfonlayer}
\end{tikzpicture}
\end{center}
\caption{Typical elements of different types of diagram semigroups. Undirected edges represent a pair of opposite direction arrows.} 
\label{fig:diagram_semigroup_types}
\end{figure}

For a finite set $X$ a \emph{diagram} is a subset of $(X\cup X')\times (X\cup X')$ where $|X|=|X'|=n$, the \emph{degree} of the diagram, and $X\cap X'=\varnothing$. 
Pictorially, we draw the points from $X$ on an upper row with those from $X'$ below, and we draw a directed edge $a\to b$ for each pair $(a,b)$ from the diagram.  For example, with $X=\{1,2,3,4,5\}$, the diagram
$\big\{
(2,1),
(2,3'),
(5,4'),
(5,5'),
(1',1),
(2',2'),
(2',3),
(3',4'),
(4',3'),
(5',5)
\big\}$
is pictured as the top diagram of Figure \ref{fig:diagram_semigroup_types}.

\begin{figure}[h]
\begin{center}
\begin{tikzpicture}
\tikzstyle{blackdot}=[draw=black,circle,fill=black,inner sep=1pt]
\tikzstyle{arrow}=[thick,->,>=angle 60]
\node [blackdot] at (0,0) (u1) {};
\node [blackdot,right of=u1] (u2) {};
\node [blackdot,right of=u2] (u3) {};
\node [blackdot,right of=u3] (u4) {};
\node [blackdot,below of=u1] (d1) {};
\node [blackdot,below of=u2] (d2) {};
\node [blackdot,below of=u3] (d3) {};
\node [blackdot,below of=u4] (d4) {};
\draw [arrow] (u1) edge (d3);
\draw [arrow,bend right] (u3) edge (u4);
\draw [arrow] (d1) edge (u2);
\draw [arrow] (d4) edge (u3);
\node at (-.75,-.5) {$\alpha$};
\begin{pgfonlayer}{background layer}
\fill  [gray6] plot (-.3,0) rectangle (3.3,-1);
\end{pgfonlayer}
\begin{scope}[yshift=-1cm]
\tikzstyle{arrow}=[white,thick,->,>=angle 60]
\node [blackdot] at (0,0) (u1) {};
\node [blackdot,right of=u1] (u2) {};
\node [blackdot,right of=u2] (u3) {};
\node [blackdot,right of=u3] (u4) {};
\node [blackdot,below of=u1] (d1) {};
\node [blackdot,below of=u2] (d2) {};
\node [blackdot,below of=u3] (d3) {};
\node [blackdot,below of=u4] (d4) {};
\draw [arrow] (u3) edge (d2);
\draw [arrow,bend right] (u3) edge (u4);
\draw [arrow,bend right] (d2) edge (d1);
\draw [arrow] (d2) edge (u1);
\node at (-.75,-.5) {$\beta$};
\begin{pgfonlayer}{background layer}
\fill  [gray6] plot (-.3,0) rectangle (3.3,-1);
\end{pgfonlayer}
\end{scope}
\begin{scope}[yshift=-.5cm,xshift=5cm]
\tikzstyle{arrow}=[black,thick,->,>=angle 60]
\node [blackdot] at (0,0) (u1) {};
\node [blackdot,right of=u1] (u2) {};
\node [blackdot,right of=u2] (u3) {};
\node [blackdot,right of=u3] (u4) {};
\node [blackdot,below of=u1] (d1) {};
\node [blackdot,below of=u2] (d2) {};
\node [blackdot,below of=u3] (d3) {};
\node [blackdot,below of=u4] (d4) {};
\draw [arrow] (u1) edge (d2);
\draw [arrow,bend right] (u1) edge (u3);
\draw [arrow,bend right] (u3) edge (u4);
\draw [arrow] (d2) edge (u2);
\draw [arrow,bend right] (d2) edge (d1);
\node at (-.75,-.5) {$\alpha\beta$};
\begin{pgfonlayer}{background layer}
\fill  [gray6] plot (-.3,0) rectangle (3.3,-1);
\end{pgfonlayer}
\end{scope}
\end{tikzpicture}
\end{center}
\caption{Combining partitioned binary relations $\alpha$ and $\beta$. The arrows of the product are induced by paths of the stacked diagram with the property that the consecutive arrows have alternating colors.}
\label{fig:PartBinRelMult}
\end{figure}
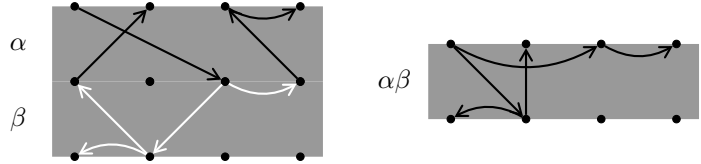

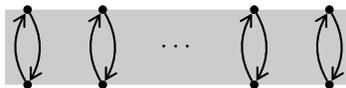
\begin{figure}
\begin{center}
\begin{center}
\begin{tikzpicture}
\tikzstyle{blackdot}=[draw=black,circle,fill=black,inner sep=1pt]
\tikzstyle{arrow}=[thick,->,>=angle 60]
\node [blackdot] at (0,0) (u1) {};
\node [blackdot,right of=u1] (u2) {};
\node [blackdot,right=1.9cm of u2] (u4) {};
\node [blackdot,right of=u4] (u5) {};
\node [blackdot,below of=u1] (d1) {};
\node [blackdot,below of=u2] (d2) {};

\node [blackdot,below of=u4] (d4) {};
\node [blackdot,below of=u5] (d5) {};

\draw [arrow,bend left] (d1) edge (u1);
\draw [arrow,bend left] (u1) edge (d1);
\draw [arrow,bend left] (d2) edge (u2);
\draw [arrow,bend left] (u2) edge (d2);
\draw [arrow,bend left] (d4) edge (u4);
\draw [arrow,bend left] (u4) edge (d4);
\draw [arrow,bend left] (d5) edge (u5);
\draw [arrow,bend left] (u5) edge (d5);
\node at (2,-.5) {$\cdots$};
\begin{pgfonlayer}{background layer}
\fill  [grey] plot (-.3,0) rectangle (4.3,-1);
\end{pgfonlayer}
\end{tikzpicture}
\end{center}
\end{center}
\caption{The identity partitioned binary relation.}
\label{fig:idpartbinrel}
\end{figure}

The product $\alpha\beta$ of two diagrams $\alpha$ and $\beta$ (on the same set $X$) is calculated as follows.  We first modify $\alpha$ and $\beta$ by changing every lower vertex $x'$ of $\alpha$ and every upper vertex $x$ of $\beta$ to $x''$.  We then stack these modified diagrams together with $\alpha$ above $\beta$ so that the vertices $x''$ are identified in the middle row (there may now be parallel edges in this stacked graph).  Finally, for each $a,b\in X\cup X'$ we include the edge $a\to b$ in $\alpha\beta$ if and only if there is a path from $a$ to $b$ in the stacked graph (as defined above) for which the edges used in the path alternate between the edges of $\alpha$ and the edges of $\beta$.  An example is given in Figure \ref{fig:PartBinRelMult} (where, for convenience, the edges of $\beta$ are white so that the kinds of paths referred to above are alternating in colour).  This operation is associative, so the set of all diagrams on the set $X$ forms a semigroup.  When $X=\{1,2,\ldots,n\}$, we denote this semigroup by $\PartBinRel_n$.  The identity element of $\PartBinRel_n$ is the diagram containing the edges $x\to x'$ and $x'\to x$ for each $x$ (see Figure \ref{fig:idpartbinrel}).

Subsemigroups of $\Partition\BinRel_n$ are the \emph{diagram semigroups} of degree $n$.
Imposing different sets of constraints on the diagrams gives rise to different kinds of diagram semigroups (example diagrams are shown in Figure~\ref{fig:diagram_semigroup_types}). There are two main ways to specialize the diagrams. We can restrict the arrows to go only one way from top to bottom, domain to codomain, yielding binary relations and then different functions, the classical transformation semigroups \cite{ClassicalTransSemigroups2009}. We can also consider partitioned binary relations that are equivalence relations, yielding the partition monoid and its submonoids.

\subsection{Binary Relations}

Prohibiting edges within the upper and lower sets and restricting to top-down edges yields $\BinRel_n$, the monoid of binary relations of an $n$-element set \cite{plemmons1970}. 

It is tempting to think that degree $n$ partitioned binary relations can be represented by binary relations of degree $2n$. This is true on the level of elements, but not on the semigroup level, since multiplication is different.
\subsection{Partial and Total Transformations}

Further constraints give us partial transformations (for each $x$, there is at most one edge $x\to y'$); transformations (for each $x$, there is exactly one edge $x\to y'$); partial permutations (injective partial transformations); and permutations (injective transformations).  The sets of all such elements are, respectively: the partial transformation semigroup $\PartialTrans_n$, the (full) transformation semigroup $\FullTrans_n$; the symmetric inverse semigroup $\InvMon_n$; and the symmetric group $\Symmetric_n$.

\subsection{Symmetric Group and Symmetric Inverse Monoid}

Degree $n$ permutations form the symmetric group $\Symmetric_n$, a central and thoroughly studied algebraic structure (e.g.~\cite{CameronPermGroups99, DixonMortimerPermGroups96}).
Partial permutations form the \emph{symmetric inverse monoid} \cite{lawson1998inverse}.

\subsection{Partitions}

If the underlying relation is an equivalence relation then we have the \emph{partition monoid} $\Partition_n$ \cite{PartitionAlgebras2005}, also known as the \emph{bipartition monoid}.
When drawing its diagrams we can omit loop edges due to reflexivity, and a pair of directed edges can be replaced by an undirected one due to the relation being symmetric.
We also use the transitivity of the equivalence relation and draw fewer edges.

\subsection{Dual Symmetric Inverse Monoid} 

The equivalence relation of a partition diagram $\alpha$ defines quotient sets of $\Dom(\alpha)$ and $\Codom(\alpha)$, sets of \emph{blocks}.
If the diagram induces a bijection between the upper and the lower set of blocks, then the diagram is a \emph{block bijection}. The monoid consisting of these block bijections is called the dual symmetric inverse monoid, since it is the categorical dual of the symmetric inverse monoid \cite{DualSymmetricInverse1998}.

\subsection{Brauer Diagrams}

Restricting to partitions of size 2 only, we get the Brauer monoid $\Brauer_n$.

\subsection{Temperley-Lieb Diagrams}

Restricting to planar diagrams from $\Brauer_n$ we get the Temperley-Lieb monoid $TL_n$.

\subsection{The big picture}

\begin{figure}
\begin{tikzpicture}
\tikzstyle{plain}=[rounded corners=3pt, draw]
\draw node [plain](Tn) {$\FullTrans_n$};
\draw node [plain,right of=Tn] (In) {$\InvMon_n$};
\draw node [plain,right of=In] (Ins) {$\DualInvMon_n$};
\draw node [plain,right of=Ins] (Br) {$\Brauer_n$};
\draw node [plain,below of=In] (Sn) {$\Symmetric_n$};

\draw node [plain,below right of=Sn] (1n) {$1_n$};
\draw node [plain,below of=Br] (TemperleyLieb) {$TL_n$};
\draw node [plain,above of=Tn] (PTn) {$\PartialTrans_n$};
\draw node [plain,above of=PTn] (Bin) {\BinRel$_n$};
\draw node [plain,above right =5mm of Bin] (PBn) {$\PartBinRel_n$};
\draw node at (3,2) [plain] (PartitionMonoid) {$\Partition_n$};

\draw  (Tn) -- (Sn);
\draw  (In) -- (Sn);
\draw  (Ins) -- (Sn);

\draw  (PartitionMonoid) -- (In);
\draw  (PartitionMonoid) -- (Ins);
\draw  (PartitionMonoid) -- (Tn);
\draw  (PartitionMonoid) -- (Br);

\draw  (TemperleyLieb) -- (1n);
\draw  (Br) -- (TemperleyLieb);
\draw  (Sn) -- (1n);
\draw  (Sn) -- (Br);

\draw (PTn) -- (Tn);
\draw (PTn) -- (In);
\draw (Bin) -- (PTn);
\draw (PBn) -- (Bin);
\draw (PBn) -- (PartitionMonoid);
\end{tikzpicture}
\begin{tabular}{cl}
\hline
Symbol & Name and References\\ 
\hline 
$\PartBinRel_n$& Partitioned binary relations \cite{PartBinRel2013}\\
\myrowcolour 
$\BinRel_n$& Binary relations \cite{plemmons1970}\\
$\PartialTrans_n$& Partial transformation semigroup \cite{ClassicalTransSemigroups2009}\\
\myrowcolour 
$\Partition_n$&(Bi)partition monoid \cite{PartitionAlgebras2005}\\
$\Brauer_n$&Brauer monoid\\
\myrowcolour 
$\Symmetric_n$&Symmetric group\cite{CameronPermGroups99, DixonMortimerPermGroups96}\\
$\FullTrans_n$&Full transformation semigroup \cite{ClassicalTransSemigroups2009}\\
\myrowcolour 
$\InvMon_n$&Symmetric inverse monoid \cite{lawson1998inverse}\\
$\DualInvMon_n$&Dual symmetric inverse monoid \cite{DualSymmetricInverse1998}\\
\myrowcolour 
$TL_n$ & Temperley-Lieb, Jones monoid $J_n$\\
\end{tabular}
\caption{Summary table of diagram semigroups: notation, names and references. \todo{refs for all}}
\label{fig:names}
\end{figure}
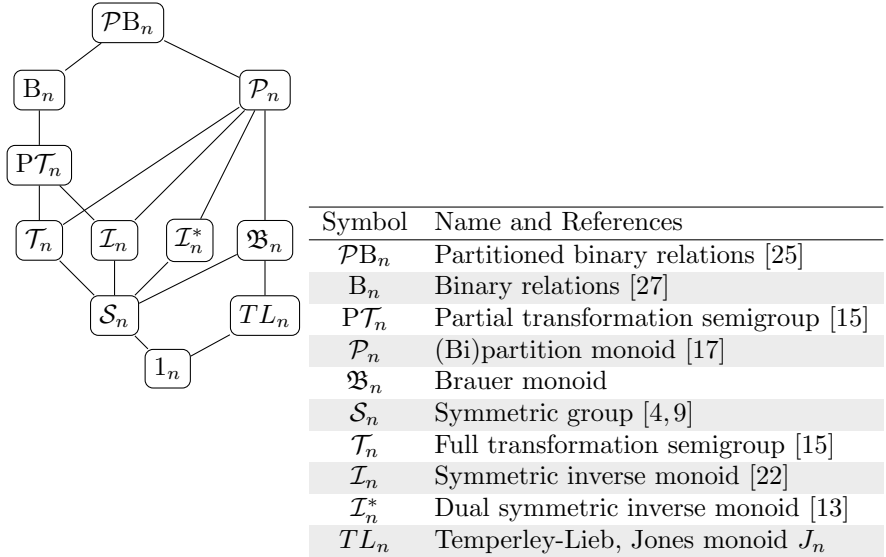

Transformations and partial permutations and the dual symmetric inverse monoid $\DualInvMon_n$ also embed into $\Partition_n$, but $P\FullTrans_n$ does not (its elements can be realized with diagrams of $\Partition_n$, but the multiplication is different).
Relationships between the diagram monoids are shown by the Hasse diagram in Figure~\ref{fig:names}, with edges denoting restrictions of rules of definition and consequently embeddings.

Some  results about diagram representations of semigroups are known.
For instance, it is trivial that the semigroup of \emph{partial transformations} $\PartialTrans_n$ embeds into the semigroup of transformations $\FullTrans_{n+1}$, so we can avoid partial maps by adding one more point.
Less trivially, it has recently been shown that in the weaker sense of generating the corresponding pseudovarieties, the Brauer monoid $\Brauer_n$ can represent all monoids, while the Temperley-Lieb monoid $\TemperleyLieb_n$ can represent all aperiodic semigroups \cite{BrauerTypeMonoids}. 

The size of each kind of monoid grows quickly with the degree for diagram semigroups (Tab.\ \ref{tab:sequences}), so brute-force enumeration of all subsemigroups is only possible for small semigroups.

\begin{table}
\begin{tabular}{clrrrrrr}
\hline
 &Order &$n=1$ & $n=2$ & $n=3$ & $n=4$ & $n=5$ & $n=6$ \\ 
\hline 
$\PartBinRel_n$& $2^{(2n)^2}$ & $16$ & 65536 &$2^{36}$ & $2^{64}$ & $2^{100}$& $2^{144}$ \\
\myrowcolour 
$\BinRel_n$& $2^{n^2}$ & 2 & 16&512&65536&$2^{25}$& $2^{36}$ \\
$\Partition_n$&$B_{2n}=\sum_1^{2n} S(2n, k)$ & 2 & 15&  203& 4140&115975 & 4213597\\
\myrowcolour 
$\PartialTrans_n$& $(n+1)^n$ & 2& 9& 64& 625& 7776&117649\\
$\DualInvMon_n$&$\sum_1^n k!\big(S(n, k)\big)^2$& 1& 3& 25& 339& 6721& 179643\\
\myrowcolour 
$\FullTrans_n$&$n^n$ & 1&4&27&256&3125& 46656\\
$\InvMon_n$&$\sum_0^n k!{n\choose k}^2$ &2& 7& 34& 209& 1546& 13327\\
\myrowcolour 
$\Brauer_n$&$(2n-1)!!$ & 1& 3& 15& 105& 945&10395\\
$\Symmetric_n$&$n!$ & 1& 2& 6& 24& 120& 720\\
\myrowcolour 
$TL_n,J_n$ & $C_n=\frac{1}{n+1} {2n \choose n}$ & 1&2&5&14&42& 132\\
\hline
\end{tabular}
\caption{Summary table of diagram semigroups. ($C_n$ is the $n$-th Catalan number, $S(n,k)$ are the Stirling numbers of the second kind, $B_n$ is the $n$-th Bell number)}
\label{tab:sequences}
\end{table}

\section{Semigroup Enumeration}
The enumeration of semigroups by computers started very early in computing history, and continuing efforts were focused on constructing abstract semigroups by finding all associative multiplication tables up to isomorphism and anti-isomorphism  of the given size \cite{For55, tamura1, tamura2, Ple67,KRS76,JW77,SZT94,monoidenum2009,ord10semigroups}.
The generated data sets are conveniently available  in the \GAP~\cite{GAP4} package called \smallsemi~\cite{smallsemi}.

As a next step, following the success story of computational group theory, where permutation group representations have efficient algorithms, basic algorithms for calculating with finite transformation semigroups were developed.  These cover, for instance, multiplying transformations, enumerating elements, deciding membership, and calculating the divisibility relations and the principal ideals, the so called $\cD$-class structure~\cite{Linton1998aa}.
Currently \GAP~\cite{GAP4} and its \semigroups~package \cite{Semigroups} have the implementations of these algorithms.

In algebraic automata theory, finite automata are represented as transformation semigroups. 
Interest in studying finite computations led to the enumeration of transformation semigroups up to degree 4 by enumerating all subsemigroups of $\FullTrans_4$.
This was achieved by dividing up the semigroup along its ideal structure, allowing parallel processing \cite{SubSemi2014}.
An \emph{ideal} is a subsemigroup $I\leq S$ such that $SI\subseteq I$ and $IS\subseteq I$.
By using Rees-quotients~\cite{clifford_preston,Howie95}, we can collapse an ideal to a zero element, substantially reducing the search space by separating it into two parts, the ideal $I$ and the quotient semigroup $S/I$.
This algorithm is implemented in the \subsemi~package \cite{subsemi}. 

The strategy of taking the all enveloping full structure and enumerating its substructures can be used for all kinds of diagram representations.
Moreover, recently the fundamental semigroup algorithms have been generalized to partial permutation semigroups, partition monoids, matrix semigroups, and subsemigroups of finite regular Rees matrix and
$0$-matrix semigroups~\cite{Semigroups,citrus}. These two facts made the computational enumeration of finite diagram semigroups possible.   

The simplest way of compressing enumeration data is to consider only conjugacy class representatives. Two diagrams are conjugate if they differ only by a reordering of their points. We denote the set of conjugacy class representatives of the semigroup $S$ by $\Sub_G(S)$, where $G$ is the permutation group of all permutations of the underlying points that preserves the semigroup $S$. 

\section{Visualising the Database}
\begin{table}
\begin{tabular}{rrrrrrr}
\hline
& $n=1$ & $n=2$ & $n=3$ & $n=4$ & $n=5$ & $n=6$\\
\hline
$\PartBinRel_n$ & 1262 & & & & & \\
\myrowcolour
$\BinRel_n$ & 4 & 385 &  &  & & \\
$\Partition_n$ & 4 & 272 &  &  & & \\
\myrowcolour 
$\PartialTrans_n$ & 4 & 50 & 94232 & & & \\
$\InvMon_n$ & 4 & 23 & 2963 & & & \\
\myrowcolour 
$\DualInvMon_n$ & 2 & 6 & 795 & & & \\
$\FullTrans_n$ & 2 & 8 & 283 &132069776 & & \\
\myrowcolour 
$\Brauer_n$ & 2 & 6 & 42 & 10411 & & \\
$\TemperleyLieb_n$ & 2 &4 &12 &232& 12592 & 324835618\\ 
\myrowcolour
$\Symmetric_n$ & 1 &2 &4 &11& 19 & 56\\ 
\end{tabular}
\caption{Summary table of the numbers of distinct (up conjugacy) semigroups of given diagram type and degree.}
\label{tab:limits}
\end{table}

The first thing we would like to know about a type of combinatorial structure is to know how many there are.
The known values for the numbers of diagram semigroups of degree $n$  are summarized in Table \ref{tab:limits}.
These values can be used for testing other methods of calculating these numbers and for devising a closed formula, if ever possible.
However, a single number does not tell us much about the type of semigroups.

\subsection{Size Distributions}
\begin{figure}
\begin{center}
\input{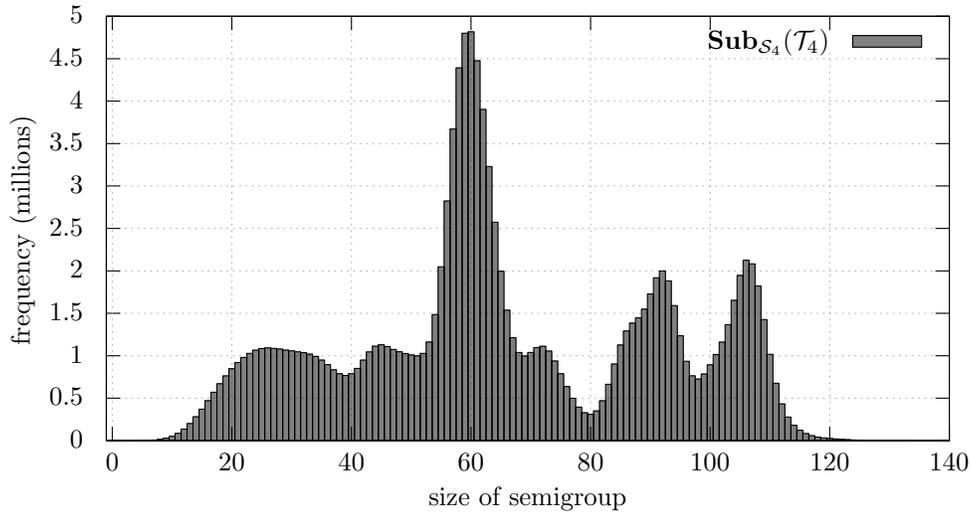}
\end{center}
\caption{The six peaks of $\FullTrans_4$. The main bulk of the size distribution of transformation semigroups of degree 4.} 
\label{fig:T4size6peaks}
\end{figure}
It is a basic fact of group theory (Lagrange's Theorem) that the subgroups of a group $G$ have orders that are divisors of $|G|$.
Therefore, the size of a permutation group of degree $n$ should be a divisor of $n!$, yielding a size distribution of these permutation groups with spikes at these values (how many groups for a divisor is a nontrivial matter).
For transformation semigroups we do not have this restriction, but still, the continuous looking curve for the size distribution of $\Sub_{\Symmetric_4}(\FullTrans_4)$ is somewhat surprising, see Fig.~\ref{fig:T4size6peaks}.

\begin{figure}
\begin{center}
\input{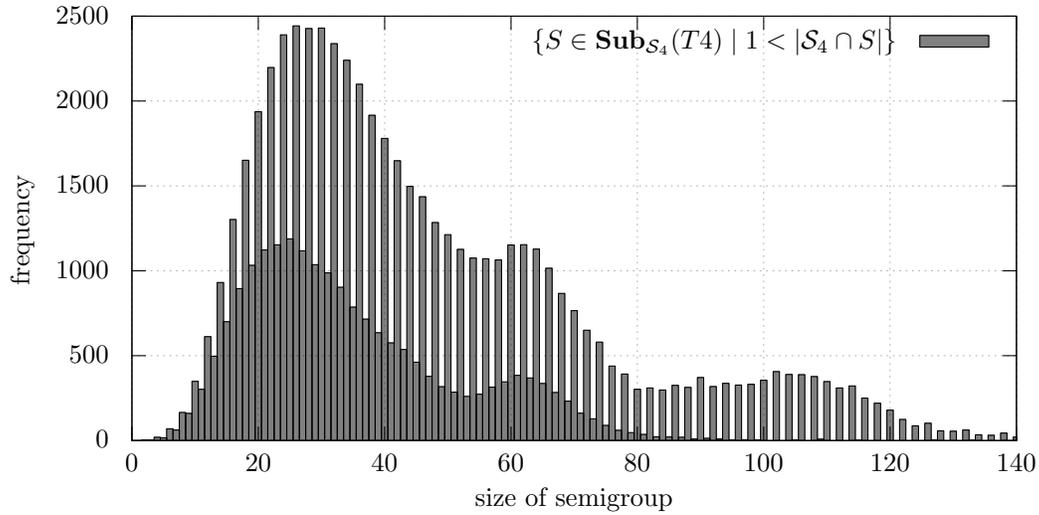}
\end{center}
\caption{The double distribution of transformation semigroups of degree 4 containing nontrivial permutations. Semigroups of even order are more abundant.}
\label{fig:PT4sizes}
\end{figure}

The next question is whether the size distribution is more `group-like' for semigroups containing nontrivial permutations.
Figure \ref{fig:PT4sizes} shows that the curve is again continuous looking. 
Moreover, it has a curious even-odd alternating pattern, actually drawing two curves. 
The pattern is mainly due to the semigroups containing only one nontrivial permutation (93.93\% of semigroups). 


\subsection{Semigroup Structure Summary Heatmaps}

When studying a semigroup, it is a standard first step to ask about its structure in terms of its Green's equivalence relations \cite{clifford_preston,Howie95,QBook}. These can be defined 
 in terms of divisibility relations
$$ t\ \cR\ s \iff \exists p,q\in S^\id \text{ such that } t=sp \text{ and } s=tq,$$
$$ t\ \cL\ s \iff \exists p,q\in S^\id \text{ such that } t=ps \text{ and } s=qt,$$
$$ t\ \cJ\ s \iff \exists p,q,u,v\in S^\id \text{ such that } t=psq \text{ and } s=utv.$$
The $\mathcal D$ relation is defined to be the composition of the $\mathcal L$ and $\mathcal R$ relations (in either order); so $t \mathcal D s$ if and only if $t\mathcal L u \mathcal R s$ for some $u\in S$.  For finite semigroups, the $\mathcal D$ relation coincides with the $\mathcal J$ relation.  Since semigroup elements are $\mathcal D$-related precisely when they can be obtained from each other by multiplication within the semigroup, the $\mathcal D$-classes can be thought of as the `local pools of reversibility'.

\subsubsection{Size versus the number of $\cD$-classes}

\begin{figure}
\input{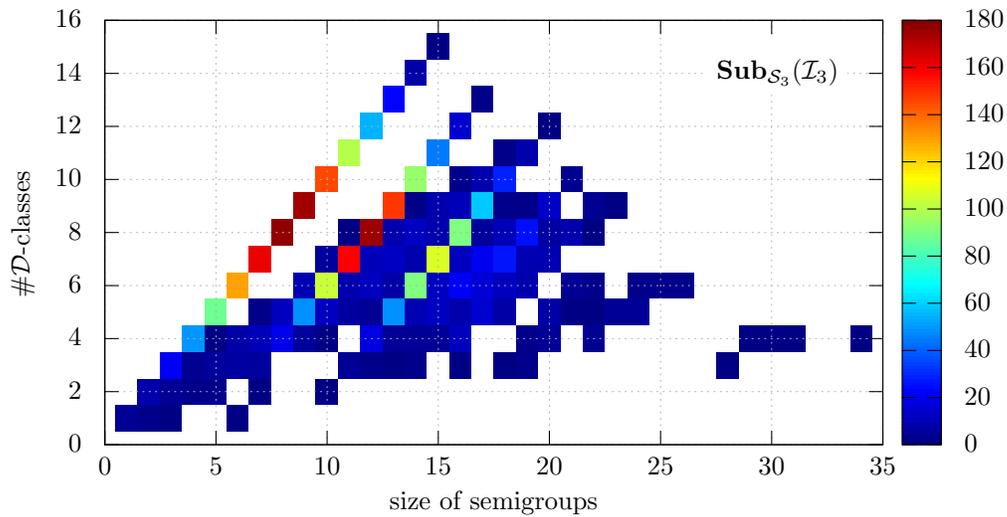}
\caption{Size versus the number of $\cD$-classes for partial permutation semigroups up to degree 3.}
\label{fig:I3SvsD}
\end{figure}

\begin{figure}
\input{tikz/T4_SvsD.tikz}
\caption{Size versus the number of $\cD$-classes for transformation semigroups up to degree 4. Frequency values in millions.}
\label{fig:T4SvsD}
\end{figure}

Each $\cD$-class of a finite semigroup contains elements that are mutually reachable from each other by left or right multiplication. 
Is there a relation between the number of elements and the number of $\cD$-classes beyond the trivial constraint that the former is an upper bound for the latter?
Visualising the relationship as a heatmap of all degree 3 inverse semigroups (Fig.~\ref{fig:I3SvsD}) we see that for many of them the $\cD$-classes are singletons. Also, the heatmap is discontinuous, which is probably explained by the fact that inverse semigroups, consisting of partial permutations,  are the closest to groups.
In contrast, the analogous heatmap for $\Sub_{\Symmetric_4}(\FullTrans_4)$ (Fig.~\ref{fig:T4SvsD}) looks in a sense continuous (up to about 155).

\subsubsection{Size versus the number of idempotents}

\begin{figure}
\input{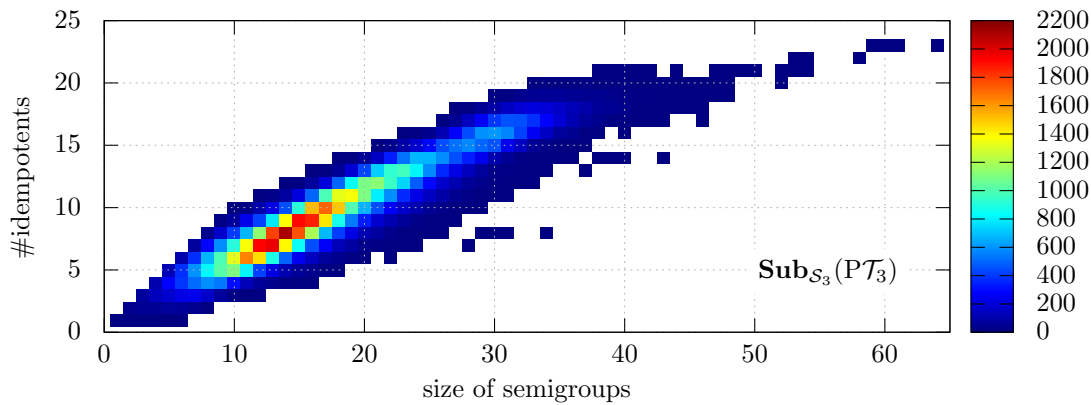}
\caption{Size versus the number of idempotents for partial transformation semigroups up to degree 3.}
\label{fig:PT3SvsI}
\end{figure}

\begin{figure}
\input{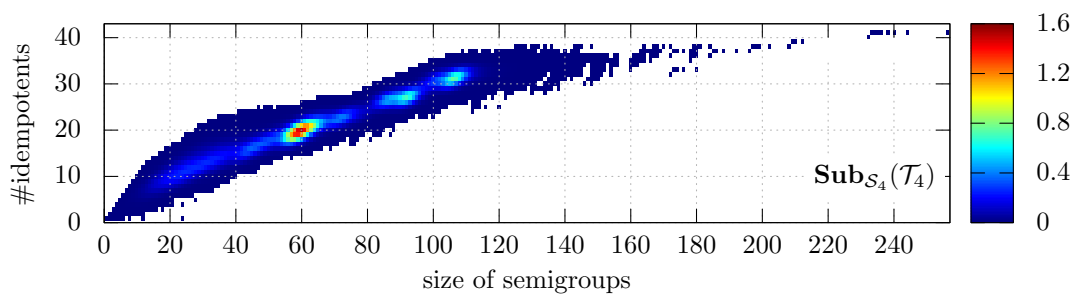}
\caption{Size versus the number of idempotents for transformation semigroups up to degree 4. Frequency values in millions.}
\label{fig:T4SvsI}
\end{figure}

Idempotents, which satisfy $x=x^2$, also play an important role in the analysis of semigroup structure since they are the identities of the subgroups of a semigroup. Figures \ref{fig:PT3SvsI} and \ref{fig:T4SvsI} show connection between the size of the semigroup and the number of idempotents. The heatmaps indicate that most of the semigroups are clustered on a line.

\subsection{Superfractals in the Temperley-Lieb Monoid}

\begin{figure}
\begin{center}
\includegraphics[width=.75\textwidth]{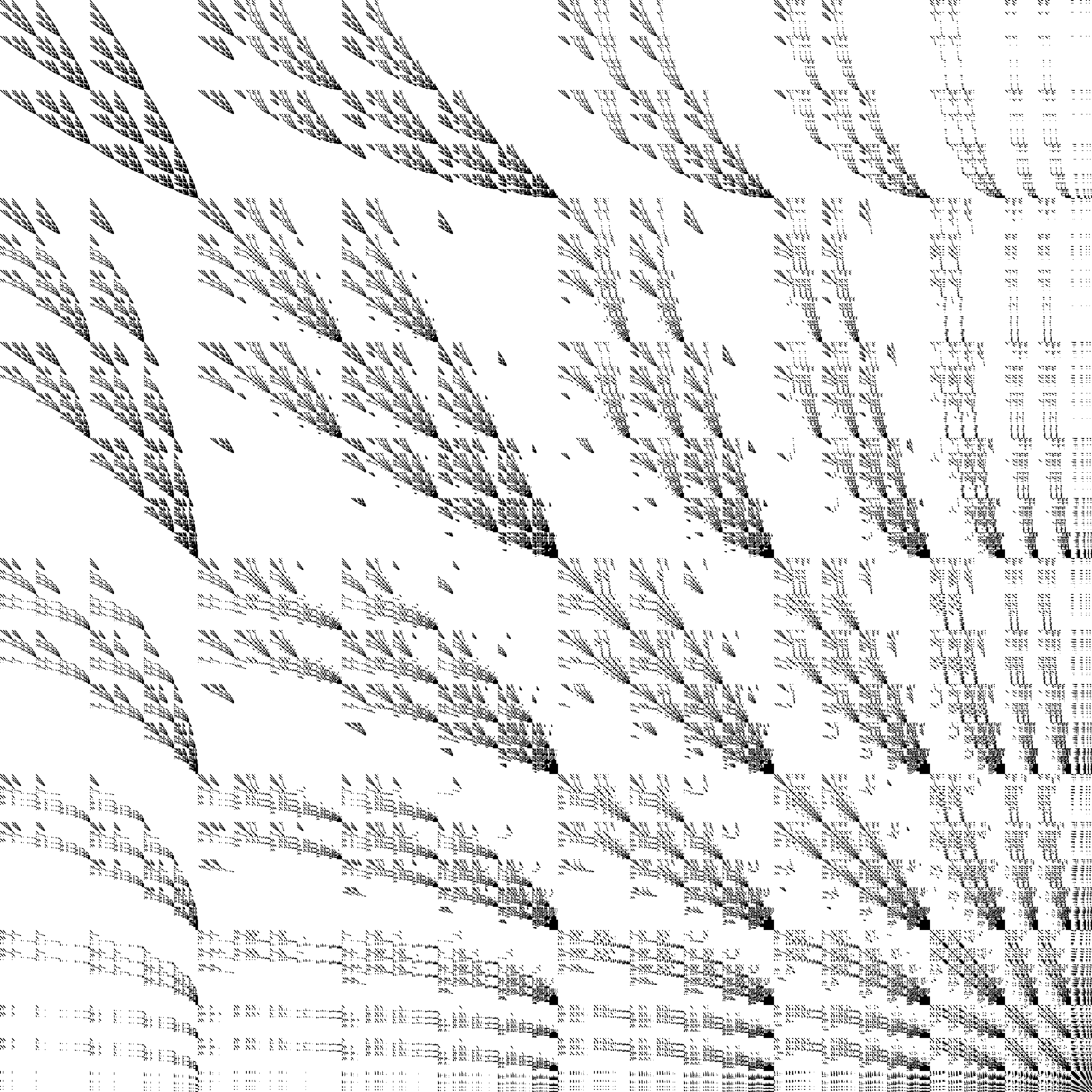}
\end{center}
\caption{The positions  of idempotents in the sixth $\cD$-class of $\TemperleyLieb_{16}$. }
\label{fig:ferns}
\end{figure}

 Looking into the structure of the individual $\cD$-classes we often draw the so called `eggbox' picture 
\cite{clifford_preston,Howie95}. 
The Temperley-Lieb monoid has a single hierarchy of $\cD$-classes determined by the number of cups and caps in the planar diagram.
However, the locations of the idempotents have an interesting structure (assuming the standard generating set), see Figure \ref{fig:ferns}. The picture looks like a superfractal \cite{barnsley2006superfractals}, but we do not know the generating iterated function systems.

\section{Open Problems}
The generated data sets are awash with interesting observations and open questions (including some not explicitly referred above). Here we just mention a few. 
\begin{enumerate}
\item Why are there six peaks in the distribution of sizes in $\Sub_{\Symmetric_4}(\FullTrans_4)$? (see Fig.~\ref{fig:T4size6peaks}) Why are semigroups of size 60 the most abundant? Is there a number theoretical or algebraic explanation? What is the shape for $\FullTrans_5$ and beyond?
\item What is the asymptotic behaviour of the ratio of subsemigroups  and  subsets of the diagram semigroups?
\item Is there a convergence in the shape of the distributions? Do they become single-peaked as the degree increases? If so, that would imply that we can talk about typical members of degree $n$ diagram semigroups in a statistical sense.
\item Explain the fractal structure appearing in $\TemperleyLieb_n$.
\end{enumerate}

\bibliography{../compsemi}

\begin{thebibliography}{10}

\bibitem{BrauerTypeMonoids}
Karl Auinger.
\newblock Pseudovarieties generated by {B}rauer type monoids.
\newblock {\em Forum Math.}, 26(1):1--24, 2014.

\bibitem{barnsley2006superfractals}
M.F. Barnsley.
\newblock {\em SuperFractals}.
\newblock Cambridge University Press, 2006.

\bibitem{Brauer1937}
Richard Brauer.
\newblock On algebras which are connected with the semisimple continuous
  groups.
\newblock {\em Ann. of Math. (2)}, 38(4):857--872, 1937.

\bibitem{CameronPermGroups99}
Peter~J. Cameron.
\newblock {\em {Permutation Groups}}.
\newblock London Mathematical Society, 1999.

\bibitem{clifford_preston}
A.H. Clifford and G.B. Preston.
\newblock {\em The Algebraic Theory of Semigroups, Vol.~1}.
\newblock Number~7 in Mathematical Surveys. American Mathematical Society, 2nd
  edition, 1967.

\bibitem{ord10semigroups}
Andreas Distler, Chris Jefferson, Tom Kelsey, and Lars Kotthoff.
\newblock The semigroups of order 10.
\newblock In Michela Milano, editor, {\em Principles and Practice of Constraint
  Programming}, Lecture Notes in Computer Science, pages 883--899. Springer
  Berlin Heidelberg, 2012.

\bibitem{monoidenum2009}
Andreas Distler and Tom Kelsey.
\newblock The monoids of orders eight, nine {\&} ten.
\newblock {\em Ann. Math. Artif. Intell.}, 56(1):3--21, 2009.

\bibitem{smallsemi}
Andreas Distler and James~D. Mitchell.
\newblock {\em {\textsc{Smallsemi} --- a data library of semigroups of small
  size, version 0.6.7}}, 2013.
\newblock
  \href{http://tinyurl.com/jdmitchell/smallsemi/}{\url{http://tinyurl.com/jdmitchell/smallsemi/}}.

\bibitem{DixonMortimerPermGroups96}
John~D. Dixon and Brian Mortimer.
\newblock {\em {Permutation Groups}}.
\newblock {Graduate Texts in Mathematics 163}. Springer, 1996.

\bibitem{citrus}
J.~East, { A. Egri-Nagy}, { J. D. Mitchell}, and Y.~P\'eresse.
\newblock Computing with semigroups.
\newblock in preparation, 2014.

\bibitem{SubSemi2014}
James East, { Attila} { Egri-Nagy}, and { James D. Mitchell}.
\newblock On enumerating transformation semigroups.
\newblock \href{http://arxiv.org/abs/1403.0274}{arXiv:1403.0274 [math.GR]},
  2014.

\bibitem{subsemi}
James East, { Attila} { Egri-Nagy}, and { James D. Mitchell}.
\newblock {\em \textsc{{S}ub{S}emi} -- software package for enumerating
  subsemigroups, Version 1.0+}, 2014.
\newblock
  \href{https://bitbucket.org/egri-nagy/subsemi}{\url{https://bitbucket.org/egri-nagy/subsemi}}.

\bibitem{DualSymmetricInverse1998}
D.~G. Fitzgerald and Jonathan Leech.
\newblock Dual symmetric inverse monoids and representation theory.
\newblock {\em Journal of the Australian Mathematical Society (Series A)},
  64:345--367, 6 1998.

\bibitem{For55}
George~E. Forsythe.
\newblock S{WAC} computes {$126$} distinct semigroups of order {$4$}.
\newblock {\em Proc. Amer. Math. Soc.}, 6:443--447, 1955.

\bibitem{ClassicalTransSemigroups2009}
Olexandr Ganyushkin and Volodymyr Mazorchuk.
\newblock {\em {Classical Transformation Semigroups}}.
\newblock {Algebra and Applications}. Springer, 2009.

\bibitem{GAP4}
The GAP~Group.
\newblock {\em {GAP -- Groups, Algorithms, and Programming, Version 4.7.7}},
  2015.
\newblock \href{http://www.gap-system.org}{\url{http://www.gap-system.org}}.

\bibitem{PartitionAlgebras2005}
T.~Halverson and A.~Ram.
\newblock Partition {A}lgebras.
\newblock {\em European J. Combin.}, 26(6):869--921, 2005.

\bibitem{Howie95}
John~M. Howie.
\newblock {\em {Fundamentals of Semigroup Theory}}, volume~12 of {\em {London
  Mathematical Society Monographs New Series}}.
\newblock Oxford University Press, 1995.

\bibitem{Jones1994}
V.~F.~R. Jones.
\newblock A quotient of the affine {H}ecke algebra in the {B}rauer algebra.
\newblock {\em Enseign. Math. (2)}, 40(3-4):313--344, 1994.

\bibitem{JW77}
H.~J{\"u}rgensen and P.~Wick.
\newblock Die {H}albgruppen der {O}rdnungen {$\leq 7$}.
\newblock {\em Semigroup Forum}, 14(1):69--79, 1977.

\bibitem{KRS76}
Daniel~J. Kleitman, Bruce~R. Rothschild, and Joel~H. Spencer.
\newblock The number of semigroups of order {$n$}.
\newblock {\em Proc. Amer. Math. Soc.}, 55(1):227--232, 1976.

\bibitem{lawson1998inverse}
M.V. Lawson.
\newblock {\em Inverse Semigroups: The Theory of Partial Symmetries}.
\newblock World Scientific, 1998.

\bibitem{Linton1998aa}
S.~A. Linton, G.~Pfeiffer, E.~F. Robertson, and N.~Ru{\v{s}}kuc.
\newblock Groups and actions in transformation semigroups.
\newblock {\em Math. Z.}, 228(3):435--450, 1998.

\bibitem{Martin1994}
Paul Martin.
\newblock Temperley-{L}ieb algebras for nonplanar statistical mechanics---the
  partition algebra construction.
\newblock {\em J. Knot Theory Ramifications}, 3(1):51--82, 1994.

\bibitem{PartBinRel2013}
Paul Martin and Volodymyr Mazorchuk.
\newblock Partitioned binary relations.
\newblock {\em Mathematica Scandinavica}, 113(1):30--52, 2013.

\bibitem{Semigroups}
James Mitchell.
\newblock {\em Semigroups Version 2.2}, 2015.
\newblock
  \href{http://www-groups.mcs.st-andrews.ac.uk/~jamesm/semigroups.php}{\url{http://www-groups.mcs.st-andrews.ac.uk/~jamesm/semigroups.php}}.

\bibitem{plemmons1970}
R.~J. Plemmons and M.~T. West.
\newblock On the semigroup of binary relations.
\newblock {\em Pacific Journal of Mathematics}, 35(3):743--753, 1970.

\bibitem{Ple67}
Robert~J. Plemmons.
\newblock There are {${\rm 15973}$} semigroups of order {${\rm 6}$}.
\newblock {\em Math. Algorithms}, 2:2--17, 1967.

\bibitem{QBook}
John Rhodes and Benjamin Steinberg.
\newblock {\em {The q-theory of Finite Semigroups}}.
\newblock Springer, 2008.

\bibitem{SZT94}
S.~Satoh, K.~Yama, and M.~Tokizawa.
\newblock Semigroups of order {$8$}.
\newblock {\em Semigroup Forum}, 49(1):7--29, 1994.

\bibitem{tamura2}
Takayuki Tamura.
\newblock Notes on finite semigroups and determination of semigroups of order
  {$4$}.
\newblock {\em J. Gakugei. Tokushima Univ. Math.}, 5:17--27, 1954.

\bibitem{tamura1}
Kazutoshi Tetsuya, Takao Hashimoto, Tadao Akazawa, Ryoichi Shibata, Tadashi
  Inui, and Takayuki Tamura.
\newblock All semigroups of order at most {$5$}.
\newblock {\em J. Gakugei Tokushima Univ. Nat. Sci. Math.}, 6:19--39. Errata on
  loose, unpaginated sheet, 1955.

\end{thebibliography}
\bibliographystyle{plain}
\end{document}